\documentclass[11pt]{amsart}
\usepackage{amsfonts}

\usepackage{graphicx}
\usepackage{amscd}

\newtheorem{theorem}{Theorem}

\newtheorem{corollary}[theorem]{Corollary}

\newtheorem{definition}{Definition}
\newtheorem{example}{Example}

\newtheorem{proposition}{Proposition}

\begin{document}
\title{On mappings of terms determined by hypersubstitutions}
\author{J\"{o}rg Koppitz}
\address{J\"{o}rg Koppitz \\
University of Potsdam\\
Institute of Mathematics\\
Postfach 601553\\
14415 Potsdam, Germany}
\email{koppitz@rz.uni-potsdam.de}
\author{Slavcho Shtrakov}
\address{Slavcho Shtrakov\\
South-West-University Blagoevgrad\\
Faculty of Mathematics and Natural Sciences\\
2700 Blagoevgrad, Bulgaria}
\email{shtrakov@aix.swu.bg.}
\date{}

\begin{abstract}
The extensions of hypersubstitutions are mappings on the set of all terms.
In the present paper we characterize all hypersubstitutions which provide
bijections on the set of all terms. The set of all such hypersubstitutions
forms a monoid.\par
On the other hand, one can modify each hypersubstitution to any mapping on
the set of terms. For this we can consider mappings $\rho $ from the set of
all hypersubstitutions into the set of all mappings on the set of all terms.
If for each hypersubstitution $\sigma $ the application of $\rho (\sigma )$
to any identity \ in a given variety $V$ is again an identity in $V$, so
that variety is called $\rho $-solid. The concept of a $\rho $-solid variety
generalizes the concept of a solid variety. In the present paper, we
determine all $\rho $-solid varieties of semigroups for particular mappings $%
\rho $.
\end{abstract}

\maketitle

\section{Basic Definitions and Notations}

We fix a type $\tau =(n_{i})_{i\in I}$, $n_{i}>0$ for all $i\in I$, and a
set of operation \ symbols $\Omega :=\{f_{i}\mid i\in I\}$ where $f_{i}$ is $n_{i}$-ary. 
Let $W_{\tau }(X)$ be the set of all terms of type $\tau $ over
some fixed alphabet $X=\{x_{1},x_{2},\ldots \}$. Terms in $W_{\tau }(X_{n})$
with $X_{n}=\{x_{1},\ldots ,x_{n}\}$, $n\geq 1$, are called $n$-ary. For
natural numbers $m,$ $n\geq 1$ we define a mapping $S_{m}^{n}:W_{\tau
}(X_{n})\times W_{\tau }(X_{m})^{n}\rightarrow W_{\tau}(X_{m})$ in the
following way: For $(t_{1},\ldots ,t_{n})\in W_{\tau }(X_{m})^{n}$ we put:%
\\
$
\begin{array}{ll}
\text{(i)} & S_{m}^{n}(x_{i},t_{1},\ldots ,t_{n}):=t_{i}\text{ for }1\leq
i\leq n;
\\
\text{(ii)} & S_{m}^{n}(f_{i}(s_{1},\ldots ,s_{n_{i}}),t_{1},\ldots
,t_{n}):=f_{i}(S_{m}^{n}(s_{1},t_{1},\ldots ,t_{n}),\ldots
,S_{m}^{n}(s_{n_{i}},t_{1}, \\
& \ldots ,t_{n}))\text{ for }i\in I,s_{1},\ldots ,s_{n_{i}}\in W_{\tau
}(X_{n})\text{ where }S_{m}^{n}(s_{1},t_{1},\ldots ,t_{n}), \\
& \ldots, S_{m}^{n}(s_{n_{i}},t_{1},\ldots ,t_{n})\text{ will be assumed to be
already defined.}
\end{array}
$\par
If it is obvious what $m$ is, we write $S^{n}$. For $t\in W_{\tau }(X)$ we
define the depth of $t$ in the following inductive way:\\
$\begin{array}{ll}
\text{(i)} & depth(t):=0\text{ for }t\in X;
\\
\text{(ii)} & depth(t):=\max \{depth(t_{1}),\ldots ,depth(t_{n_{i}})\}+1\\ &
\text{ for }t=f_{i}(t_{1},\ldots ,t_{n_{i}}) \text{with }i\in I, t_{1},\ldots ,t_{n_{i}}\in W_{\tau }(X) \text{
where } \\
& depth(t_{1}),\ldots ,depth(t_{n_{i}})\text{ will  }  \text{be assumed to be already defined}.\par
\end{array}
$\par
By $c(t)$ we denote the length of a term $t$ (i.e. the number of the
variables occurring in $t)$, $var(t)$ denotes the set of all variables
occurring in $t$ and $cv(t)$ means the number of elements in the set $var(t)$.
 Instead of $x_{1},x_{2},x_{3},\ldots $ we write also $x,y,z,\ldots $.

The concept of a hypersubstitution was introduced in \cite{d2}.

\begin{definition}
A mapping $\sigma :\Omega \rightarrow W_{\tau }(X)$ which assigns to every $%
n_{i}$-ary operation symbol $f_{i}$, $i\in I$, an $n_{i}$-ary term is called
a hypersubstitution of type $\tau $ (shortly hypersubstitution).
The set of all hypersubstitutions of type $\tau $ will be denoted by $%
Hyp(\tau )$.
\end{definition}

To each hypersubstitution $\sigma $ there belongs a \ mapping from the set
of all terms of the form $f_{i}(x_{1},\ldots ,x_{n_{i}})$ to the terms $%
\sigma (f_{i})$. It follows that every hypersubstitution of type $\tau $
then induces a mapping $\widehat{\sigma }:W_{\tau }(X)\rightarrow W_{\tau
}(X)$ as follows:\\
$
\begin{array}{ll}
\text{(i)} & \widehat{\sigma }[w]:=w\text{ for }w\in X;
\\
\text{(ii)} & \widehat{\sigma }[f_{i}(t_{1},\ldots ,t_{n_{i}})%
]:=S^{n}(\sigma (f_{i}),\widehat{\sigma }[t_{1}],\ldots ,\widehat{\sigma }[%
t_{n_{i}}])\text{ for }i\in I,t_{1},\ldots ,t_{n_{i}}\text{ } \\
& \in W_{\tau }(X)\text{ where }\widehat{\sigma }[t_{1}],\ldots ,\widehat{%
\sigma }[t_{n_{i}}]\text{ will be assumed to be already } \\ & \text{  defined.}
\end{array}
$

By $\sigma _{1}\circ _{h}\sigma _{2}$ $:=\widehat{\sigma }_{1}\circ \sigma
_{2}$ is defined an associative operation on $Hyp(\tau )$ where $\circ $
denotes the usual composition of mappings. By $\varepsilon $ we denote the
hypersubstitution with $\varepsilon (f_{i})=f_{i}(x_{1},\ldots ,x_{n_{i}})$
for $i\in I$, where $\varepsilon $ deals as identity element. Then $%
(Hyp(\tau );\circ _{h},\varepsilon )$ forms a monoid, denoted by $\mathbf{%
Hyp(\tau )}$.

\section{Bijections on $W_{\protect\tau }(X)$}

By $Bij(\tau )$ we denote the set of all $\sigma \in Hyp(\tau )$ such that $%
\widehat{\sigma }:W_{\tau }(X)\rightarrow W_{\tau }(X)$ is a bijection on $%
W_{\tau }(X)$. Such hypersubstitutions have a high importance in computer
science.

The product of two bijections is again a bijection. Further, for two
hypersubstitutions $\sigma _{1}$ and $\sigma _{2}$ we have
\begin{equation*}
(\sigma _{1}\circ _{h}\sigma _{2})\widehat{}=\widehat{\sigma }_{1}\circ
\widehat{\sigma }_{2}
\end{equation*}
(see \cite{d3}). So we have the following result.

\begin{proposition}
($Bij(\tau );\circ _{h},\varepsilon )$ forms a submonoid of $\mathbf{%
Hyp(\tau )}$.
\end{proposition}

For the characterization of $Bij(\tau )$ we need the following notations:%

(i) $\mathcal{B}$ denotes the set of all bijections on $\Omega $ preserving
the arity.

(ii) Let $S_{n}$ be the set of all permutations of the set $\{1,\ldots ,n\}$\\
for $1\leq n\in \mathbb{N}$.

(iii) $A:=\underset{1\leq n\in \mathbb{N}}{\bigcup }S_{n}.$ \par
(iv) $\mathcal{P}:=\{p\in A^{I}\mid p(i)\in S_{n_{i}}$ for $i\in I\}.$

The following theorem characterizes $Bij(\tau )$ for any type $\tau $.

\begin{theorem}
Let $\tau =(n_{i})_{i\in I}$, $n_{i}>0$ for all $i\in I$, be any type. For
each $\sigma \in Hyp(\tau )$ the following statements are equivalent:\par
$
\begin{array}{ll}
\text{(i)} & \sigma \in Bij(\tau ).
\\
\text{(ii)} & \text{There are }h\in \mathcal{B}\text{ and }p\in \mathcal{P}%
\text{ such that} \\
& \sigma (f_{i})=h(f_{i})(x_{p(i)(1)},\ldots ,x_{p(i)(n_{i})})\text{ for all
}i\in I\text{.}
\end{array}
$
\end{theorem}

\begin{proof}
$(ii)\Rightarrow (i):$ We show by induction that $\widehat{\sigma }$ is
injective and surjective. \par
Injectivity: Let $s,t\in W_{\tau }(X)$ with $\widehat{\sigma }[s]=\widehat{%
\sigma }[t]$. \par
Suppose that the $depth(s)=0.$ Then $depth(t)=0$ and $s,t$ are variables
with $s=$ $\widehat{\sigma }[s]=\widehat{\sigma }[t]=t$. \par
Suppose that from $\widehat{\sigma }[s`]=\widehat{\sigma }[t`]$ there
follows $s`=t`$ for any $s`,t`\in W_{\tau }(X)$ with $depth(s`)\leq n$.%
\par
Let $depth(s)=n+1$. Then $depth(t)\geq 1$ and there are $i,j\in I$ with $%
s=f_{i}(s_{1},\ldots ,s_{n_{i}})$ and $t=f_{j}(t_{1},\ldots ,t_{n_{j}})$.
Now we have \par
$\widehat{\sigma }[s]=S^{n_{i}}(h(f_{i})(x_{p(i)(1)},\ldots
,x_{p(i)(n_{i})}),\widehat{\sigma }[s_{1}],\ldots ,\widehat{\sigma }[%
s_{n_{i}}])$ and \par
$\widehat{\sigma }[t]=S^{n_{j}}(h(f_{j})(x_{p(j)(1)},\ldots
,x_{p(j)(n_{j})}),\widehat{\sigma }[t_{1}],\ldots ,\widehat{\sigma }[%
t_{n_{j}}])$. From $\widehat{\sigma }[s]=\widehat{\sigma }[t]$ it follows
that $h(f_{i})=h(f_{j})$ and thus $f_{i}=f_{j}$, i.e. $i=j$, since $h$ is a
bijection. Hence $S^{n_{i}}(h(f_{i})(x_{p(i)(1)},\ldots ,x_{p(i)(n_{i})}),%
\widehat{\sigma }[s_{1}],\ldots ,\widehat{\sigma }[s_{n_{i}}])$\par
$=S^{n_{i}}(h(f_{i})(x_{p(i)(1)},\ldots ,x_{p(i)(n_{i})}),\widehat{\sigma }[%
t_{1}],\ldots ,\widehat{\sigma }[t_{n_{j}}])$ and,\\  consequently, $\widehat{%
\sigma }[s_{k}]=\widehat{\sigma }[t_{k}]$ for $1\leq k\leq n_{i}$. By our
hypothesis we get $s_{k}=t_{k}$ for $1\leq k\leq n_{i}$. Consequently, $%
s=f_{i}(s_{1},\ldots ,s_{n_{i}})=f_{j}(t_{1},\ldots ,t_{n_{j}})=t$.\par
Surjectivity: For $w\in X$ we have $\widehat{\sigma }[w]=w$. \par
Suppose that for any $s\in W_{\tau }(X)$ with $depth(s)\leq n$ there is an $%
\widetilde{s}\in W_{\tau }(X)$ with $\widehat{\sigma }[\widetilde{s}]=s$.
\par
Let now $t\in W_{\tau }(X)$ be a term with $depth(t)=n+1$. Then there is an $%
i\in I$ with $t=f_{i}(t_{1},\ldots ,t_{n_{i}})$ and by our hypothesis there
are $\widetilde{t}_{1},\ldots ,\widetilde{t}_{n_{i}}\in W_{\tau }(X)$ such
that $\widehat{\sigma }[\widetilde{t}_{k}]=t_{k}$ for $1\leq k\leq n_{i}$.
Further there is a $j\in I$ with $h(f_{j})=f_{i}$ and $n_{i}=n_{j}$. Now we
consider the term $\widetilde{t}:=f_{j}(\widetilde{t}_{p(j)^{-1}(1)},\ldots ,%
\widetilde{t}_{p(j)^{-1}(n_{i})})$. There holds $$\widehat{\sigma }[%
\widetilde{t}]=S^{n_{i}}(h(f_{j})(x_{p(j)(1)},\ldots ,x_{p(j)(n_{j})}),%
\widehat{\sigma }[\widetilde{t}_{p(j)^{-1}(1)}],\ldots ,\widehat{\sigma }[%
\widetilde{t}_{p(j)^{-1}(n_{i})}])$$
$$=S^{n_{i}}(f_{i}(x_{p(j)(1)},\ldots
,x_{p(j)(n_{i})}),t_{p(j)^{-1}(1)},\ldots ,t_{p(j)^{-1}(n_{i})})$$ (by
hypothesis) 
$$=f_{i}(t_{1},\ldots ,t_{n_{i}})=t.$$
\par
$(i)\Rightarrow (ii):$ Since $\widehat{\sigma }$ is surjective for each $%
j\in I$ there is an $s_{j}\in W_{\tau }(X)$ with $\widehat{\sigma }[s_{j}%
]=f_{j}(x_{1},\ldots ,x_{n_{j}})$ which is minimal with respect to the
depth. Obviously, the case $\ depth(s_{j})=0$ is impossible. Thus there are
a $k\in I$ and $r_{1},\ldots ,r_{n_{k}}\in W_{\tau }(X)$ with $%
s_{j}=f_{k}(r_{1},\ldots ,r_{n_{k}})$. So $$\widehat{\sigma }[s_{j}]=\widehat{%
\sigma }[f_{k}(r_{1},\ldots ,r_{n_{k}})]=S^{n_{k}}(\sigma (f_{k}),\widehat{%
\sigma }[r_{1}],\ldots ,\widehat{\sigma }[r_{n_{k}}])
=f_{j}(x_{1},\ldots ,x_{n_{j}}).$$ This is only possible if $\sigma
(f_{k})\in X$ or $\sigma (f_{k})=f_{j}(a_{1},\ldots ,a_{n_{j}})$ with $$%
a_{1},\ldots ,a_{n_{j}}\in \{x_{1},\ldots ,x_{n_{k}}\}, \mid
\{a_{1},\ldots ,a_{n_{j}}\}\mid =n_{j},$$ and thus $n_{k}\geq n_{j}$. But the
case $\sigma (f_{k})\in X$ is impossible. Otherwise there is an $i\in
\{1,\ldots ,n_{k}\}$ with $\sigma (f_{k})=x_{i}$ and $\widehat{\sigma }[s_{j}%
]=\widehat{\sigma }[r_{i}]$ where $depth(s_{j})>depth(r_{i})$, this
contradicts the minimallity of $s_{j}$. This shows that for all $j\in I$
there are a $k(j)\in I$ with $n_{k(j)}\geq n_{j}$ and $a_{1},\ldots
,a_{n_{j}}\in X_{n_{k(j)}}$ with $\mid \{a_{1},\ldots ,a_{n_{j}}\}\mid
=n_{j} $ such that $\sigma (f_{k(j)})=f_{j}(a_{1},\ldots ,a_{n_{j}})$.%
\par
Assume that $n_{k(j)}>n_{j}$ for some $j\in I$. Then there is an $x\in
X_{n_{k(j)}}\setminus var(\sigma (f_{k(j)}))$, i.e. $x$ is not essential in $%
\sigma (f_{k(j)})$ and thus $\widehat{\sigma }$ is not a bijection on $%
W_{\tau }(X)$ (see [1], [6]), a contradiction. Thus $n_{k(j)}=n_{j}$ and $%
\sigma (f_{k(j)})=f_{j}(x_{\pi _{j}(1)},\ldots ,x_{\pi _{j}(n_{j})})$ for
some $\pi _{j}\in S_{n_{j}}$. \par
Assume that there are $j,l\in I$ with $l\neq k(j)$ such that $f_{j\text{ }}$%
is the first operation symbol in $\sigma (f_{l})$. We put $t:=\widehat{%
\sigma }[f_{l}(x_{1},\ldots ,x_{n_{l}})]$. Then $t=f_{j}(t_{1},\ldots
,t_{n_{j}})$ for some $t_{1},\ldots ,t_{n_{j}}\in W_{\tau }(X)$. Since $%
\widehat{\sigma }$ is surjective, there are $s_{1},\ldots ,s_{n_{j}}\in
W_{\tau }(X)$ with $\widehat{\sigma }[s_{i}]=t_{i}$ for $1\leq i\leq n_{j}$.
Then $\widehat{\sigma }[f_{k(j)}(s_{\pi _{j}^{-1}(1)},\ldots ,s_{\pi
_{j}^{-1}(n_{j})})]$\par
$=S^{n_{j}}(\sigma (f_{k(j)}),\widehat{\sigma }[s_{\pi _{j}^{-1}(1)}],\ldots
,\widehat{\sigma }[s_{\pi _{j}^{-1}(n_{j})}])$\par
$=S^{n_{j}}(f_{j}(x_{\pi _{j}(1)},\ldots ,x_{\pi _{j}(n_{j})}),t_{\pi
_{j}^{-1}(1)},\ldots ,t_{\pi _{j}^{-1}(n_{j})})$\par
$=f_{j}(t_{1},\ldots ,t_{n_{j}})$. \par
Since $f_{k(j)}(s_{\pi _{j}^{-1}(1)},\ldots ,s_{\pi _{j}^{-1}(n_{j})})\neq
f_{l}(x_{1},\ldots ,x_{n_{l}})$, $\widehat{\sigma }$ is no injective, a
contradiction. Altogether this shows that the mapping $h:\Omega \rightarrow
\Omega $ where $h(f)$ is the first operation symbol in $\sigma (f)$ is a
bijection on $\Omega $ preserving the arity. Further, let $p\in A^{I}$ with $%
p(i):=\pi _{i}$ for $i\in I$. Then $p\in \mathcal{P}.$\par
Consequently, we have $\sigma (f_{i})=h(f_{i})(x_{p(i)(1)},\ldots
,x_{p(i)(n_{i})})$ for all $i\in I$.
\end{proof}

Let us give the following examples.

\begin{example}
Let $2\leq $ $n\in \mathbb{N}$. We consider the type $\tau _{n}=(n)$, where $%
f$ denotes the $n$-ary operation symbol. For $\pi \in S_{n}$ we define:%
\par
\begin{equation*}
\sigma _{\pi }:f\mapsto f(x_{\pi (1)},\ldots ,x_{\pi (n)})\text{.}
\end{equation*}
\par
These hypersubstitutions are precisely the bijections, i.e. $Bij(\tau
_{n})=\{\sigma _{\pi }\mid \pi \in S_{n}\}$.\par
In particular, if $n=2$ then $Bij(\tau _{2})=\{\varepsilon ,\sigma _{d}\}$
where $\sigma _{d}$ is defined by
\begin{equation*}
\sigma _{d}:f\mapsto f(x_{2},x_{1})\text{.}
\end{equation*}
\end{example}

\begin{example}
Let now $\tau =(2,2)$ where $f$ and $g$ are the both binary operation
symbols. Then we define the following eight hypersubstitutions $\sigma
_{1},\ldots ,\sigma _{8}$ by:
\begin{equation*}
\begin{tabular}{lll}
& $f\mapsto $ & $g\mapsto $ \\ \hline
$\sigma _{1}:$ & $f(x_{1},x_{2})$ & $g(x_{1},x_{2})$ \\
$\sigma _{2}:$ & $f(x_{1},x_{2})$ & $g(x_{2},x_{1})$ \\
$\sigma _{3}:$ & $f(x_{2},x_{1})$ & $g(x_{1},x_{2})$ \\
$\sigma _{4}:$ & $f(x_{2},x_{1})$ & $g(x_{2},x_{1})$ \\
$\sigma _{5}:$ & $g(x_{1},x_{2})$ & $f(x_{1},x_{2})$ \\
$\sigma _{6}:$ & $g(x_{1},x_{2})$ & $f(x_{2},x_{1})$ \\
$\sigma _{7}:$ & $g(x_{2},x_{1})$ & $f(x_{1},x_{2})$ \\
$\sigma _{8}:$ & $g(x_{2},x_{1})$ & $f(x_{2},x_{1})$%
\end{tabular}
\end{equation*}

 These hypersubstitutions are precisely the bijections, so $$Bij(\tau
)=\{\sigma _{1},\ldots ,\sigma _{8}\}.$$
\end{example}

\section{$\protect\rho $-solid varieties}

In Section 1, we mentioned that any hypersubstitution $\sigma $ can be
uniquely extended to a mapping $\widehat{\sigma }:W_{\tau }(X)\rightarrow
W_{\tau }(X)$ ($\widehat{\sigma }\in W_{\tau }(X)^{W_{\tau }(X)})$. Thus a
mapping $\rho :Hyp(\tau )\rightarrow W_{\tau }(X)^{W_{\tau }(X)}$ is defined
by setting $\rho (\sigma )=\widehat{\sigma }$ for all $\sigma \in Hyp(\tau )$%
. \par
In \cite{d4}, the concept of a solid variety was introduced. By Birkhoff, a
variety $V$ is a class of algebras of type $\tau $ satisfying a set $\Sigma $
of identities, i.e. $V=Mod\Sigma $. For a variety $V$ of type $\tau $ we
denote by $IdV$ the set of all identities in $V$. The variety $V$ is said to
be solid iff $\widehat{\sigma }[s]\approx \widehat{\sigma }[t]\in IdV$ for
all $s\approx t\in IdV$ and all $\sigma \in Hyp(\tau )$. \ For a submonoid $%
\mathbf{M}$ of $\mathbf{Hyp(\tau )}$, the variety $V$ is said to be $M$%
-solid iff $\widehat{\sigma }[s]\approx \widehat{\sigma }[t]\in IdV$ for all
$s\approx t\in IdV$ and all $\sigma \in M$ (see \cite{d3}). If $M=$ $Hyp(\tau )$
then we have solid varieties.\par
In this section we will study mappings $\rho :Hyp(\tau )\rightarrow W_{\tau
}(X)^{W_{\tau }(X)}$ and generalize the concept of an $M$-solid variety to
the concept of an $M$-$\rho $-solid variety. For convenience, we put $\sigma
^{\rho }:=\rho (\sigma )$ for $\sigma \in Hyp(\tau )$.

\begin{definition}
Let $\rho :Hyp(\tau )\rightarrow W_{\tau }(X)^{W_{\tau }(X)}$ be a mapping
and $V$ be a variety of type $\tau $ and $\mathbf{M}$ be a submonoid of $%
\mathbf{Hyp(\tau )}.$ $V$ is called $M$-$\rho $-solid iff $\sigma ^{\rho
}(s)\approx \sigma ^{\rho }(t)\in IdV$ for all $s\approx t\in IdV$ and all $%
\sigma \in M$.\par
If $M=$ $Hyp(\tau )$ then $V$ is said to be $\rho $-solid.
\end{definition}

\begin{example}
Let $\rho :Hyp(\tau )\rightarrow W_{\tau }(X)^{W_{\tau }(X)}$ be defined by $%
\rho (\sigma )=\widehat{\sigma }$ for all $\sigma \in Hyp(\tau )$. Then the $%
\rho $-solid varieties are exactly the solid varieties, which is clear by
the appropriate definitions. L. Pol\'{a}k has determined all solid varieties
of semigroups in \cite{d5}. Besides the trivial variety, exactly the self-dual
varieties in the interval between the normalization $Z\vee RB$ of the
variety of all rectangular bands and the variety defined by the identities $%
x^{2}\approx x^{4}$, $x^{2}y^{2}z\approx x^{2}yx^{2}yz$, $xy^{2}z^{2}\approx
xyz^{2}yz^{2}$, and $xyzyx\approx xyxzxyx$ as well as the varieties $RB$ of
all rectangular, $NB$ of all normal, and $RegB$ of all regular bands are
solid.
\end{example}

\bigskip

In Section 2 we have checked that $Bij(\tau )$ forms a monoid. For
particular mappings $\rho :Hyp(\tau )\rightarrow W_{\tau }(X)^{W_{\tau }(X)}$
the $Bij(\tau )$-$\rho $-solid varieties are of special interest, in
particular for type $\tau =(2)$ and semigroup varieties. They realize
substitutions of operations in terms which are useful in some calculational
aspects of computer algebra systems. In the following we will consider such
mappings $\rho :Hyp(\tau )\rightarrow W_{\tau }(X)^{W_{\tau }(X)}$.

\begin{definition}
Let
\begin{equation*}
fa:Hyp(\tau )\rightarrow W_{\tau }(X)^{W_{\tau }(X)}\text{ and }sa:Hyp(\tau
)\rightarrow W_{\tau }(X)^{W_{\tau }(X)}
\end{equation*}
be the following mappings: For $\sigma \in Hyp(\tau )$ we put\\
$
\begin{array}{ll}
\text{(i)} & \sigma ^{fa}(x):=\sigma ^{sa}(x):=x\text{ for }x\in X;\par
\\
\text{(ii)} & \sigma ^{fa}(f_{i}(t_{1},\ldots ,t_{n_{i}})):=S^{n_{i}}(\sigma
(f_{i}),\sigma ^{sa}(t_{1}),\ldots ,\sigma ^{sa}(t_{n_{i}}))\text{ and} \\
& \sigma ^{sa}(f_{i}(t_{1},\ldots ,t_{n_{i}})):=f_{i}(\sigma
^{fa}(t_{1}),\ldots ,\sigma ^{fa}(t_{n_{i}}))\text{ for }i\in I\text{ and}
\\
& t_{1},\ldots ,t_{n_{i}}\in W_{\tau }(X)\text{ where }\sigma
^{sa}(t_{1}),\ldots ,\sigma ^{sa}(t_{n_{i}}),\sigma ^{fa}(t_{1}),\ldots
,\sigma ^{fa}(t_{n_{i}}) \\
& \text{will be assumed to be already defined.}
\end{array}
$
\end{definition}

If we consider $M$-$\rho $-solid varieties of semigroups we have the type $%
\tau =(2)$ and thus $\rho :Hyp(2)\rightarrow W_{(2)}(X)^{W_{(2)}(X)}$ (where
$Hyp(2):=Hyp((2))$). If one considers semigroup identities, we have the
associative law and we can renounce of the operation symbol $f$ and the
brackets, i.e. we write semigroup words only as sequences of variables.

\begin{theorem}
The trivial variety $TR$ and the variety $Z$ of all zero semigroups (defined
by $xy\approx zt$) are the only $sa$-solid varieties of semigroups.
\end{theorem}

\begin{proof}
Clearly, $TR$ is $sa$-solid.\par
We show that for any $\sigma \in Hyp(2)$ and any $t\in W_{(2)}(X)$ there
holds $\sigma ^{sa}(t)\approx t\in IdZ.$\par
If $t\in X$ then $\sigma ^{sa}(t)=t$.\par
If $t\notin X$ then $t=f(t_{1},t_{2})$ for some $t_{1},t_{2}\in $ $%
W_{(2)}(X) $. Thus $c(t)\geq 2$ and $t\approx xy\in IdZ.$ Further, there
holds $\sigma ^{sa}(t)=f(\sigma ^{fa}(t_{1}),\sigma ^{fa}(t_{2}))\approx
xy\in IdZ$. Consequently, $\sigma ^{sa}(t)\approx t\in IdZ$. \par
This shows that $\sigma ^{sa}(s)\approx s\approx t\approx \sigma ^{sa}(t)$
holds in $Z$ for all $s\approx t\in IdZ$ and all $\sigma \in Hyp(2)$, i.e. $%
Z $ is $sa$-solid.\par
\par
Conversely, let $V$ be an $sa$-solid variety of semigroups. By $\sigma _{x}$
($\sigma _{y}$) we will denote the hypersubstitution which maps the binary
operation symbol $f$ to the term $x_{1}$ ($x_{2}$). Then $\sigma
_{x}^{sa}(f(f(x,y),z))\approx \sigma _{x}^{sa}(f(x,f(y,z)))\in IdV$. This
provides $xz\approx xy\in IdV$. From $\sigma _{y}^{sa}(f(f(x,y),z))\approx
\sigma _{y}^{sa}(f(x,f(y,z)))\in IdV$ it follows $yz\approx xz\in IdV$. Both
identities $xz\approx xy$ and $yz\approx xz$ provide $yz\approx xt$, i.e. $%
V\subseteq Z$. But $TR$ and $Z$ are the only subvarieties of $Z$.
\end{proof}

\begin{proposition}
A variety $V$ of semigroups is $Bij(2)$-$sa$-solid iff \par
(i) $V\subseteq Mod\{x(yz)\approx (xy)z,\ xyz\approx zxy\}$ and \par
(ii) $V\subseteq Mod\{x(yz)\approx (xy)z,\ xyz\approx xzy\approx zxy\}$ if
there is an identity $s\approx t\in IdV$ with $cv(s)=c(s)=3$ and $c(t)\neq 3$
or $cv(t)\neq 3$ or $var(t)\neq var(s)$.
\end{proposition}

\begin{proof}
We have already mentioned that $Bij(2)=\{\varepsilon ,\sigma _{d}\}$.
\par
Suppose that $V$ is $Bij(2)$-$sa$-solid. Then $$\sigma
_{d}^{sa}(f(f(x,y),z))\approx \sigma _{d}^{sa}(f(x,f(y,z)))\in
IdV,$$ so $yxz\approx xzy\in IdV$. Let now $s\approx t\in IdV$ with $%
cv(s)=c(s)=3$. \par
If $c(t)\leq 2$ then $\sigma _{d}^{sa}(t)=t.$\par
If $c(t)\geq 4$ then $\sigma _{d}^{sa}(t)\approx t\in IdV$ is easy to check
using $yxz\approx xzy\in IdV$.\par
If $c(t)=3$ and $cv(t)=1$ then $\sigma _{d}^{sa}(t)\approx t\in IdV$ is
obvious. \par
If $c(t)=3$ and $cv(t)=2$ then there are $w_{1},w_{2}\in X$ such that $%
t=(w_{1}w_{2})w_{2}$ or \ $t=(w_{2}w_{1})w_{2}$ or $t=(w_{2}w_{2})w_{1}$ or $%
t=w_{1}(w_{2}w_{2})$ or $t=w_{2}(w_{1}w_{2})$ or $t=w_{2}(w_{2}w_{1})$.
Using $yxz\approx xzy\in IdV$ we get that $w_{1}w_{2}w_{2}\approx
w_{2}w_{1}w_{2}\approx w_{2}w_{2}w_{1}$ in $V$. This shows that $\sigma
_{d}^{sa}(t)\approx t\in IdV$.\par
From $cv(s)=c(s)=3$ it follows $s=(w_{1}w_{2})w_{3}$ or $\
s=w_{1}(w_{2}w_{3})$ for some $w_{1},w_{2},w_{3}\in X$. Without loss of
generality let $s=w_{1}(w_{2}w_{3})$, so $\sigma
_{d}^{sa}(s)=w_{1}w_{3}w_{2} $.\par
If $c(t)\neq 3$ or $cv(t)\neq 3$, from $\sigma _{d}^{sa}(s)\approx \sigma
_{d}^{sa}(t)\in IdV$ it follows $w_{1}w_{3}w_{2}\approx t\in IdV$.
Consequently, $w_{1}w_{3}w_{2}\approx w_{1}w_{2}w_{3}\in IdV$. \par
If $cv(t)=c(t)=3$ and $var(t)\neq var(s)$ then there is a $w\in
var(t)\setminus var(s)$. Substituting $w$ by $w^{2}$ we get $s\approx r\in
IdV$ from $s\approx t\in IdV$ where $c(r)=4$. Then we get $xyz\approx zxy\in
IdV$ as above.\par
\par
Suppose that (i) and (ii) are satisfied. Let $s\approx t\in IdV$. Then $%
\varepsilon ^{sa}(s)\approx \varepsilon ^{sa}(t)\in IdV$. We have to show
that $\sigma _{d}^{sa}(s)\approx \sigma _{d}^{sa}(t)\in IdV$ and consider
the following cases:\par
1) If $c(s)\neq 3$ or $cv(s)\neq 3$ and $c(t)\neq 3$ or $cv(t)\neq 3$ then
we have $\sigma _{d}^{sa}(s)\approx s\in IdV$ and $\sigma
_{d}^{sa}(t)\approx t\in IdV$ as we have shown already. This provides $%
\sigma _{d}^{sa}(s)\approx \sigma _{d}^{sa}(t)\in IdV$. \par
2.1) If $cv(s)=c(s)=3$ and $c(t)\neq 3$ or $cv(t)\neq 3$ or $var(t)\neq
var(s)$ then $xyz\approx xzy\approx zxy$ holds in $V$ (by (ii)) and it is
easy to see that $\sigma _{d}^{sa}(s)\approx s\in IdV$ and $\sigma
_{d}^{sa}(t)\approx t\in IdV$, so $\sigma _{d}^{sa}(s)\approx \sigma
_{d}^{sa}(t)\in IdV$. \par
2.2) If $cv(s)=c(s)=3$ and $c(t)=3$ and $cv(t)=3$ and $var(t)=var(s)$ then
there are $w_{1},w_{2},w_{3}\in X$ such that $s,t\in \{r_{1},\ldots
,r_{12}\} $ where \par
\begin{tabular}{llll}
$r_{1}=w_{2}(w_{1}w_{3})$ & $r_{2}=(w_{2}w_{1})w_{3}$ & $%
r_{3}=w_{3}(w_{2}w_{1})$ & $r_{4}=(w_{3}w_{2})w_{1}$ \\
$r_{5}=w_{1}(w_{3}w_{2})$ & $r_{6}=(w_{1}w_{3})w_{2}$ & $%
r_{7}=w_{2}(w_{3}w_{1})$ & $r_{8}=(w_{2}w_{3})w_{1}$ \\
$r_{9}=w_{3}(w_{1}w_{2})$ & $r_{10}=(w_{3}w_{1})w_{2}$ & $%
r_{11}=w_{1}(w_{2}w_{3})$ & $r_{12}=(w_{1}w_{2})w_{3}.$%
\end{tabular}
\par
Then $\sigma _{d}^{sa}(r_{1})=r_{7},$ $\sigma _{d}^{sa}(r_{2})=r_{12},$ $%
\sigma _{d}^{sa}(r_{3})=r_{9},$ $\sigma _{d}^{sa}(r_{4})=r_{8},$ $\sigma
_{d}^{sa}(r_{5})=r_{11},$ $\sigma _{d}^{sa}(r_{6})=r_{10},$ $\sigma
_{d}^{sa}(r_{7})=r_{1},$ $\sigma _{d}^{sa}(r_{8})=r_{4},$ $\sigma
_{d}^{sa}(r_{9})=r_{3},$ $\sigma _{d}^{sa}(r_{10})=r_{6},$ $\sigma
_{d}^{sa}(r_{11})=r_{5},$ and $\sigma _{d}^{sa}(r_{12})=r_{2}.$ This shows
that $\sigma _{d}^{sa}(r_{i})\approx $ $\sigma _{d}^{sa}(r_{j})\in IdV$ for $%
1\leq i,j\leq 6$ or $7\leq i,j\leq 12$ by $xyz\approx zxy\in IdV$. If $%
r_{i}\approx r_{j}\in IdV$ with $1\leq i\leq 6$ or $7\leq j\leq 12$ or
conversely, then $xyz\approx xzy\in IdV$. Together with $xyz\approx zxy\in
IdV$ it is easy to check that then $\sigma _{d}^{sa}(r_{i})\approx r_{i}$ $%
\in IdV$ and $\sigma _{d}^{sa}(r_{j})\approx r_{j}\in IdV$, i.e. $\sigma
_{d}^{sa}(r_{i})\approx $ $\sigma _{d}^{sa}(r_{j})\in IdV$. Altogether this
shows that $\sigma _{d}^{sa}(s)\approx \sigma _{d}^{sa}(t)\in IdV$. \par
3) If $cv(t)=c(t)=3$ then we get dually $\sigma _{d}^{sa}(s)\approx \sigma
_{d}^{sa}(t)\in IdV$.
\end{proof}

\begin{theorem}
$TR$ is the only $fa$-solid variety of semigroups.
\end{theorem}

\begin{proof}
Clearly, $TR$ is $fa$-solid. Let $V$ be an $fa$-solid variety of semigroups.
From $\sigma _{x}^{fa}(f(f(x,y),z))\approx \sigma _{x}^{fa}(f(x,f(y,z)))\in
IdV$ it follows $xy\approx x\in IdV$. Moreover, $\sigma
_{y}^{fa}(f(f(x,y),z))\approx \sigma _{y}^{fa}(f(x,f(y,z)))\in IdV$ provides
$z\approx yz\in IdV$. Both identities $xy\approx x$ and $z\approx yz$ give $%
z\approx y$, i.e. $V=TR$.
\end{proof}

\begin{proposition}
A variety $V$ of semigroups is $Bij(2)$-$fa$-solid iff \par
(i) $V\subseteq Mod\{x(yz)\approx (xy)z,\ xyz\approx zxy\}$ and \par
(ii) $V$ is a variety of commutative semigroups if there is an identity $%
s\approx t\in IdV$ with $cv(s)=c(s)=2$ and $c(t)\neq 2$ or $cv(t)\neq 2$ or $%
var(t)\neq var(s)$.
\end{proposition}

\begin{proof}
We have already mentioned that $Bij(2)=\{\varepsilon ,\sigma _{d}\}$.
\par
Suppose that $V$ is $Bij(2)$-$fa$-solid. Then $$\sigma
_{d}^{fa}(f(f(x,y),z))\approx \sigma _{d}^{fa}(f(x,f(y,z)))\in
IdV,$$ so $zxy\approx yzx\in IdV$. Let now $s\approx t\in IdV$ with $%
cv(s)=c(s)=2$. \par
If $c(t)=1$ then $\sigma _{d}^{fa}(t)=t.$\par
If $c(t)\geq 3$ then $\sigma _{d}^{fa}(t)\approx t\in IdV$ is easy to check
using $zxy\approx yzx\in IdV$.\par
If $c(t)=2$ and $cv(t)=1$ then $\sigma _{d}^{fa}(t)\approx t\in IdV$ is
obvious.\par
From $cv(s)=c(s)=2$ it follows $s=w_{1}w_{2}$, so $\sigma
_{d}^{fa}(s)=w_{2}w_{1}$.\par
If $c(t)\neq 2$ or $cv(t)\neq 2$ from $\sigma _{d}^{fa}(s)\approx \sigma
_{d}^{fa}(t)\in IdV$ it follows $w_{2}w_{1}\approx t\in IdV$ and,
consequently, $w_{1}w_{2}\approx w_{2}w_{1}\in IdV$. \par
If $cv(t)=c(t)=2$ and $var(t)\neq var(s)$ then there is a $w\in
var(t)\setminus var(s)$. Substituting $w$ by $w^{2}$ we get $s\approx r\in
IdV$ from $s\approx t\in IdV$ where $c(r)=3$. Then we get $xy\approx yx\in
IdV$ as above.\par
\par
Suppose that (i) and (ii) are satisfied. Let $s\approx t\in IdV$. Then $%
\varepsilon ^{fa}(s)\approx \varepsilon ^{fa}(t)\in IdV$. We have to show
that $\sigma _{d}^{fa}(s)\approx \sigma _{d}^{fa}(t)\in IdV$ and consider
the following cases:\par
1) If $c(s)\neq 2$ or $cv(s)\neq 2$ and $c(t)\neq 2$ or $cv(t)\neq 2$ then
we have $\sigma _{d}^{fa}(s)\approx s\in IdV$ and $\sigma
_{d}^{fa}(t)\approx t\in IdV$ as we have shown already. This provides $%
\sigma _{d}^{fa}(s)\approx \sigma _{d}^{fa}(t)\in IdV$. \par
2.1) If $cv(s)=c(s)=2$ and $c(t)\neq 2$ or $cv(t)\neq 2$ or $var(t)\neq
var(s)$ then $V$ is a variety of commutative semigroups (by (ii)) and it is
easy to see that $\sigma _{d}^{fa}(s)\approx s\in IdV$ and $\sigma
_{d}^{fa}(t)\approx t\in IdV$, so $\sigma _{d}^{fa}(s)\approx \sigma
_{d}^{fa}(t)\in IdV$. \par
2.2) If $cv(s)=c(s)=2$ and $c(t)=$ $cv(t)=2$ and $var(t)=var(s)$ then there
are $w_{1},w_{2}\in X$ such that $s=w_{1}w_{2}$ or $s=w_{2}w_{1}$ and $%
t=w_{1}w_{2}$ or $t=w_{2}w_{1}$.\par
If $s=t$ then $\sigma _{d}^{fa}(s)=\sigma _{d}^{fa}(t).$\par
If $s\neq t$ then $s\approx t$ is the commutative law and we have $\sigma
_{d}^{fa}(s)\approx \sigma _{d}^{fa}(t)\in IdV$. \par
3) \ If $cv(t)=c(t)=2$ then we get dually $\sigma _{d}^{fa}(s)\approx \sigma
_{d}^{fa}(t)\in IdV$.
\end{proof}

\begin{definition}
We define a mapping $\gamma _{n}:Hyp(\tau )\rightarrow W_{\tau }(X)^{W_{\tau
}(X)}$ for each natural number $n$ as follows: For $\sigma \in Hyp(\tau )$
we put\\
$
\begin{array}{ll}
\text{(i)} & \sigma ^{\gamma _{0}}:=\widehat{\sigma };\par
\\
\text{(ii)} & \sigma ^{\gamma _{n}}(x):=x\text{ for }x\in X\text{ and }1\leq
n\in \mathbb{N}; \\
\text{(iii)} & \sigma ^{\gamma _{n}}(f_{i}(t_{1},\ldots
,t_{n_{i}})):=f_{i}(\sigma ^{\gamma _{n-1}}(t_{1}),\ldots ,\sigma ^{\gamma
_{n-1}}(t_{n_{i}}))\text{ for }1\leq n\in \mathbb{N}\text{,} \\
& i\in I\text{, and }t_{1},\ldots ,t_{n_{i}}\in W_{\tau }(X)\text{.}
\end{array}
$\par
We put $Hyp^{(n)}(\tau ):=\{\sigma ^{\gamma _{n}}\mid \sigma \in Hyp(\tau
)\} $ for $n\in \mathbb{N}$.
\end{definition}

For the hypersubstitution $\varepsilon \in Hyp(\tau )$ (the identity element
in $Hyp(\tau )$) there holds $\varepsilon ^{\gamma _{n}}=\widehat{%
\varepsilon }$ for all $n\in \mathbb{N}$. This becomes clear by the
following considerations: We have $\varepsilon ^{\gamma _{0}}=\widehat{%
\varepsilon }$ and suppose that $\varepsilon ^{\gamma _{n}}=\widehat{%
\varepsilon }$ for some natural number $n$ then there holds $\varepsilon
^{\gamma _{n+1}}(x)=x=\widehat{\varepsilon }[x]$ and $\varepsilon ^{\gamma
_{n+1}}(f_{i}(t_{1},\ldots ,t_{n_{i}}))$\par
$=f_{i}(\varepsilon ^{\gamma _{n}}(t_{1}),\ldots ,\varepsilon ^{\gamma
_{n}}(t_{n_{i}}))$\par
$=f_{i}(\widehat{\varepsilon }[t_{1}],\ldots ,\widehat{\varepsilon }[%
t_{n_{i}}])$\par
$=f_{i}(t_{1},\ldots ,t_{n_{i}})$\par
$=\widehat{\varepsilon }[f_{i}(t_{1},\ldots ,t_{n_{i}})]$.

\begin{proposition}
The monoids $(Hyp^{(n)}(\tau );\circ ,\widehat{\varepsilon })$ and $\mathbf{%
Hyp(\tau )}$ are isomorphic for each natural number $n.$
\end{proposition}

\begin{proof}
Let $n$ be a natural number. We define a mapping $h:Hyp(\tau )\rightarrow
Hyp^{(n)}(\tau )$ by $h(\sigma ):=\sigma ^{\gamma _{n}}$ for $\sigma \in
Hyp(\tau )$. We show that $h$ is injective. For this let $\sigma _{1},\sigma
_{2}\in Hyp(\tau )$ with $\sigma _{1}^{\gamma _{n}}=\sigma _{2}^{\gamma
_{n}} $. Assume that $\sigma _{1}\neq \sigma _{2}$. Then there is an $i\in I$
with $\sigma _{1}(f_{i})\neq \sigma _{2}(f_{i})$ and we have $\widehat{%
\sigma }_{1}[f_{i}(x_{1},\ldots ,x_{n_{i}})]\neq \widehat{\sigma }%
_{2}[f_{i}(x_{1},\ldots ,x_{n_{i}})]$. Then we define:\par
$
\begin{array}{ll}
\text{(i)} & t_{0}:=f_{i}(x_{1},\ldots ,x_{n_{i}}); \\
\text{(ii)} & t_{p+1}:=f_{i}(t_{p},x_{2},\ldots ,x_{n_{i}})\text{ for }p\in
\mathbb{N}.
\end{array}
$\par
It is easy to check that $\sigma _{1}^{\gamma _{n}}(t_{n})\neq \sigma
_{2}^{\gamma _{n}}(t_{n})$ because of $\widehat{\sigma }_{1}[t_{0}]\neq
\widehat{\sigma }_{2}[t_{0}]$, which contradicts $\sigma _{1}^{\gamma
_{n}}=\sigma _{2}^{\gamma _{n}}$. This shows that $h$ is injective.\par
Clearly, $h$ is surjective. Consequently, $h$ is a bijective mapping.
\par
It is left to show that $h$ satisfies the homomorphic property. We will show
by induction on $n$ that $h(\sigma _{1}\circ _{h}\sigma _{2})=h(\sigma
_{1})\circ h(\sigma _{2})$, i.e. $(\sigma _{1}\circ _{h}\sigma _{2})^{\gamma
_{n}\text{ }}=\sigma _{1}^{\gamma _{n}}\circ \sigma _{2}^{\gamma _{n}}$.
\par
If $n=0$ then we have $\sigma _{1}^{\gamma _{0}}\circ \sigma _{2}^{\gamma
_{0}}=\widehat{\sigma }_{1}\circ \widehat{\sigma }_{2}=(\sigma _{1}\circ
_{h}\sigma _{2})\widehat{}=(\sigma _{1}\circ _{h}\sigma _{2})^{\gamma _{0}}$
(see \cite{d3}). \par
For $n=m$ we suppose that $\sigma _{1}^{\gamma _{m}}\circ \sigma
_{2}^{\gamma _{m}}=(\sigma _{1}\circ _{h}\sigma _{2})^{\gamma _{m}\text{ }}$%
. \par
Let now $n=m+1$. Obviously, we have $(\sigma _{1}^{\gamma _{m+1}}\circ
\sigma _{2}^{\gamma _{m+1}})(x)=x=(\sigma _{1}\circ _{h}\sigma _{2})^{\gamma
_{m+1}}(x)$. \par
Let $i\in I$ and $t_{1},\ldots ,t_{n_{i}}\in W_{\tau }(X)$. Then there holds
\par
$(\sigma _{1}^{\gamma _{m+1}}\circ \sigma _{2}^{\gamma
_{m+1}})(f_{i}(t_{1},\ldots ,t_{n_{i}}))=\sigma _{1}^{\gamma
_{m+1}}(f_{i}(\sigma _{2}^{\gamma _{m}}(t_{1}),\ldots ,\sigma _{2}^{\gamma
_{m}}(t_{n_{i}})))$\par
$=f_{i}((\sigma _{1}^{\gamma _{m}}\circ \sigma _{2}^{\gamma
_{m}})(t_{1}),\ldots ,(\sigma _{1}^{\gamma _{m}}\circ \sigma _{2}^{\gamma
_{m}})(t_{n_{i}}))$\par
$=f_{i}((\sigma _{1}\circ _{h}\sigma _{2})^{\gamma _{m}}(t_{1}),\ldots
,(\sigma _{1}\circ _{h}\sigma _{2})^{\gamma _{m}}(t_{n_{i}}))$ (by
hypothesis)\par
$=(\sigma _{1}\circ _{h}\sigma _{2})^{\gamma _{m+1}\text{ }%
}(f_{i}(t_{1},\ldots ,t_{n_{i}}))$. \par
Altogether, this shows that $\sigma _{1}^{\gamma _{m+1}}\circ \sigma
_{2}^{\gamma _{m+1}}=(\sigma _{1}\circ _{h}\sigma _{2})^{\gamma _{m+1}\text{
}}$.
\end{proof}

\bigskip

By definition, a variety $V$ of type $\tau $ is $M$-$\gamma _{0}$-solid iff $%
V$ is $M$-solid. The class of all solid varieties of semigroups was
determined in \cite{d5}. We will now characterize the $\gamma _{n}$-solid
varieties of semigroups for $1\leq n\in \mathbb{N}$. Here we need some else
notations. For a fixed variable $w\in X$ we put:\par
$F_{0}:=\{f(f(x,y),z)\approx f(x,f(y,z))\}$ and \par
$F_{m+1}:=\{f(s,w)\approx f(t,w)\mid s\approx t\in F_{m}\}\cup
\{f(w,s)\approx f(w,t)\mid s\approx t\in F_{m}\}$ for $m\in \mathbb{N}$.

\begin{theorem}
Let $1\leq n\in \mathbb{N}$ and $V$ be a variety of semigroups. Then $V$ is $%
\gamma _{n}$-solid iff
\begin{equation*}
x_{1}\ldots x_{n+1}\approx y_{1}\ldots y_{n+1}\in IdV\text{.}
\end{equation*}
\end{theorem}

\begin{proof}
Suppose that $V$ is $\gamma _{n}$-solid. \par
Since the associative law is satisfied in $V$ there holds $F_{n-1}\subseteq
IdV$. Since $V$ is $\gamma _{n}$-solid the application of $\sigma
_{x}^{\gamma _{n}}$ to the identities of $F_{n-1}$ gives again identities in
$V$:
\begin{equation*}
I_{1}:=\{w^{a}xzw^{b}\approx w^{a}xyw^{b}\mid a,b\in \mathbb{N},\
a+b=n-1\}\subseteq IdV\text{.}
\end{equation*}
The application of $\sigma _{y}^{\gamma _{n}}$ to the identities of $F_{n-1}$
provides
\begin{equation*}
I_{2}:=\{w^{a}yzw^{b}\approx w^{a}xzw^{b}\mid a,b\in \mathbb{N},\
a+b=n-1\}\subseteq IdV\text{.}
\end{equation*}
It is easy to check that one can derive $x_{1}\ldots x_{n+1}\approx
y_{1}\ldots y_{n+1}$ from $I_{1}\cup I_{2}$. Thus $x_{1}\ldots
x_{n+1}\approx y_{1}\ldots y_{n+1}\in IdV$.\par
\par
Suppose now that $x_{1}\ldots x_{n+1}\approx y_{1}\ldots y_{n+1}\in IdV$. We
show that for any $\sigma \in Hyp(2)$ and any $t\in W_{(2)}(X)$ there holds $%
\sigma ^{\gamma _{n}}(t)\approx t\in IdV.$\par
If $t$ contains at most $n$ operation symbols then $\sigma ^{\gamma
_{n}}(t)=t$ by definition of the mapping \ $\sigma ^{\gamma _{n}}$. \par
If $t$ contains more than $n$ operation symbols then $c(t)\geq n+1$ and $%
t\approx x_{1}\ldots x_{n+1}\in IdV$ because of $x_{1}\ldots x_{n+1}\approx
y_{1}\ldots y_{n+1}\in IdV.$ Since $t$ contains more than $n$ operation
symbols, by definition of the mapping $\sigma ^{\gamma _{n}}$, the term $%
\sigma ^{\gamma _{n}}(t)$ contains at least $n$ operation symbols and thus $%
c(\sigma ^{\gamma _{n}}(t))\geq $ $n+1$. Using $x_{1}\ldots x_{n+1}\approx
y_{1}\ldots y_{n+1}\in IdV$ we get $\sigma ^{\gamma _{n}}(t)\approx
x_{1}\ldots x_{n+1}\in IdV$. Consequently, $\sigma ^{\gamma _{n}}(t)\approx
t\in IdV$. \par
This shows that $\sigma ^{\gamma _{n}}(s)\approx s\approx t\approx \sigma
^{\gamma _{n}}(t)$ holds in $V$ for $s\approx t\in IdV$ and $\sigma \in
Hyp(2)$, i.e. $V$ is $\gamma _{n}$-solid.
\end{proof}

\begin{corollary}
$TR$ and $Z$ are the only $\gamma _{1}$-solid varieties of semigroups.
\end{corollary}

\begin{proof}
By Theorem 4, a variety $V$ of semigroups is $\gamma _{1}$-solid iff $%
x_{1}x_{2}\approx y_{1}y_{2}\in IdV,$ i.e. $V\subseteq Z$. But $TR$ and $Z$
are the only subvarieties of $Z$.
\end{proof}

\end{document}